\pdfoutput=1
\documentclass[11pt]{amsart}

\usepackage{amsmath}
\usepackage{amssymb}
\usepackage{enumitem}

\usepackage{mystyle}

\usepackage{listings}
\usepackage[enableskew,vcentermath]{youngtab}

\usepackage[bottom]{footmisc}

\usepackage{tikz}
\usepackage{tikz-cd}

\usepackage{hyperref}
\definecolor{PineGreen}{rgb}{0.0,0.47,0.44}
\definecolor{MidnightBlue}{rgb}{0.1,0.1,0.44}
\definecolor{magenta}{rgb}{1.0,0.0,1.0}
\hypersetup{
  colorlinks=true,        
  linkcolor=MidnightBlue,     
  citecolor=PineGreen,        
  filecolor=magenta,      
  urlcolor=MidnightBlue           
}

\usepackage[norelsize,ruled]{algorithm2e}
\usepackage{cite}

\usepackage{cleveref}
\crefname{thm}{theorem}{theorems}
\Crefname{thm}{Theorem}{Theorems}
\crefname{prop}{proposition}{propositions}
\Crefname{prop}{Proposition}{Propositions}
\crefname{corollary}{Corollary}{Corollaries}
\Crefname{corollary}{Corollary}{Corollaries}

\newcommand{\xqedhere}[2]{%
  \rlap{\hbox to#1{\hfil\llap{\ensuremath{#2}}}}}

\title{Computing Segre classes in arbitrary projective varieties}

\author{Corey Harris}
\address{Department of Mathematics, Florida State University, Tallahassee FL, 32306, USA}
\email{charris@math.fsu.edu}
\urladdr{\url{http://www.math.fsu.edu/~charris}}
\thanks{The author was partially supported by NSA award H98230-15-1-0027.}

\begin{document}

\begin{abstract}
  We give an algorithm for computing Segre classes of subschemes of arbitrary projective
  varieties by computing degrees of a sequence of linear projections.
  Based on the fact that Segre classes of
  projective varieties commute with intersections by general effective Cartier divisors,
  we can compile a system of linear equations which determine the coefficients for
  the Segre class pushed forward to projective space.
  The algorithm presented here comes after several others which solve the problem
  in special cases, where the ambient variety is for instance projective space;
  to our knowledge, this is the first algorithm to be able to compute Segre classes
  in projective varieties with arbitrary singularities.
\end{abstract}

\maketitle

\section{Introduction}
\label{sec:intro}

Segre classes are of fundamental importance in Fulton-MacPherson intersection theory,
and there has been a growing recent interest in their concrete computation.
For one thing, the Segre class of the singularity subscheme of a hypersurface in
$\bbP^N$ can be used to compute its Chern-Schwartz-MacPherson ($c_{SM}$) class,
and therefore also its topological Euler characteristic.
It also features heavily in enumerative problems (see \cite[\S{}9.1]{Ful98} for examples).

For the Segre class of a closed subscheme $X \overset{i}{\hookrightarrow} \bbP^N$
of projective space, Aluffi \cite{Alu03} gave the first algorithm for the computation
of $i_* \segre(X,\bbP^N) \in A_*(\bbP^N)$.
Following this, a new algorithm (call it EJP) was contributed by Eklund, Jost and Peterson \cite{EJP13,Jos13a}.
A key feature of EJP is that it could be implemented using techniques from numerical algebraic geometry.
Following these initial developments, several improvements have been made.  
First, Moe and Qviller \cite{MQ13} extended the EJP algorithm to the computation of
$j_* \segre(X,T) \in A_*(T)$ where $X \overset{j}{\hookrightarrow} T$ is a closed
embedding in a smooth projective toric variety (with some minor restrictions on $T$).
The current state-of-the-art in this direction is the algorithm of Helmer \cite{Hel15,Hel14}
which offers significant performance improvements over previous algorithms for both
$i_*\segre(X, \bbP^N)$ and $j_*\segre(X,T)$.

In a different direction, Aluffi \cite{Alu13} found a combinatorial method for computing $\segre(X,Y)$
where $Y$ is a variety (over an algebraically closed field) and $X \subset Y$ is an
\emph{r.c.~monomial subscheme}.  This places a strong restriction on the possible
subschemes that can be considered, but allows the most freedom for the ambient scheme.
However, this method requires computations of intersection numbers of divisors on $Y$
and so requires more input than just the ideals of $X$ and $Y$.

Despite the range of results just mentioned,
it is not difficult to find simple situations in which none of these methods apply.
For instance, if ${C = \Proj \big( k[x,y,z,t]/(xy-z^2) \big)}$
is the singular cone and $L$ is the subvariety cut out by the ideal $(\bar{x},\bar{z})$,
then none of the aforementioned methods can be used to compute $\segre(L,C)$,
whereas the algorithm of the present paper puts no restrictions on the ambient variety,
so it could be used here. \\

The main result of this paper is \cref{alg:main} for computing
$i_* \segre(X,Y) \in A_*(\bbP^N)$ with $Y \subset \bbP^N$ a variety and
$X \hookrightarrow Y$ an arbitrary closed subscheme.
The algorithm is based on the fact that, under suitable hypotheses,
Segre classes commute with taking hypersurface sections (\cref{cor:segreClassesCommuteWithDivisors})
and a projection technique which produces a system of linear equations in the
coefficients of $i_* \segre(X,Y)$ (\cref{prop:segredegreeformula}).

\Cref{sec:background} begins with the basic setup and assumptions for the rest of the
paper.
In \cref{subsec:backgroundPrevious} the algorithms of \cite{Alu03} and \cite{EJP13} are
reviewed and their relation to the present work is indicated.
In \cref{sec:segreClassesCommuteWithDivisors} we study the effect on Segre classes when intersecting by divisors, and \cref{sec:segreDegrees} is devoted to the proof of
\cref{prop:segredegreeformula}, which yields a linear equation satisfied by the coefficients
of the Segre class.

Finally, in \cref{sec:algorithm} the main algorithm for the computation of Segre classes
is presented.
Then in \cref{sec:Applications} several applications of the algorithm are discussed, including computations of the Chern-Schwartz-MacPherson, Chern-Mather, and polar classes of projective varieties.
We also show how to use \cref{alg:main} to compute Euclidean distance degrees and
maximum likelihood degrees in some concrete cases.
Finally, we describe some situations where the Segre class is enough to compute intersection products.
Examples of each of these are worked out in detail.

\subsection*{Acknowledgements} Thanks to Paolo Aluffi for suggesting this problem and for his all feedback during the preparation of this work.

\section{Background}
\label{sec:background}

Let $X^r \subset Y^n \subset \bbP^N$ be closed embeddings of schemes over a field $k$.
Further assume that $Y$ is a variety, or more generally a reduced subscheme of pure dimension $n$. 
Given the blowup square\\
\begin{minipage}{\linewidth}
  \centering
  \begin{tikzcd}[ampersand replacement=\&]
    \& E \arrow{r} \arrow{d}[pos=0.5, font=\small]{\pi^{\prime}} \& \tilde{Y} \ar{d}[pos=0.5,font=\small]{\pi} \\
    \& X \ar[pos=0.25,font=\small]{r} \& Y
  \end{tikzcd}
\end{minipage}
the Segre class of $X \subsetneq Y$ may be defined as
\[ \segre(X,Y) = \pi'_* \segre(E,\tilde{Y}) = \pi'_*(\frac{[E]}{1+E} \capp [\tilde{Y}]) , \]
a class in $A_*(X)$.
In general, it will not be possible to compute $s(X,Y)$, as the Chow group for $X$ may not be known,
but if we are willing to settle for the pushforward to $\bbP^N$, then the class can be written
$i_* \segre(X,Y) = s_0 [\bbP^0] + s_1 [\bbP^1] + \dots + s_r [\bbP^r] \in A_*(\bbP^N)$,
where $i: X \hookrightarrow \bbP^N$ is the inclusion.
In this paper, we exhibit an algorithm for computing these $s_i$.

\emph{Note:} the setup just presented will remain fixed throughout the paper.  Any symbols
used above will be used exactly this way from here on.

\subsection{Previous approaches}
\label{subsec:backgroundPrevious}

We now summarize some approaches that have previously been used in the computation of Segre classes.
Let $X \subset Y \subset \bbP^N$ be closed embeddings with $X$ defined by equations
$f_0,\dots,f_k$ all of degree $\deg(f_i) = d$.
We can define a rational map $\mathrm{pr} : Y \dashrightarrow \bbP^k$

\begin{minipage}{\linewidth}
  \centering
  \begin{tikzcd}[ampersand replacement=\&]
    \& \tilde{Y} \ar{dr}{\tilde{\mathrm{pr}}} \ar{dl}[swap]{\pi} \& \\
    Y \ar{rr}[dashed]{\mathrm{pr}} \& \& \bbP^k
  \end{tikzcd}
\end{minipage}
by $pr(p) = [f_0(p):\dots:f_k(p)]$ and resolve the indeterminacies by
lifting to $\tilde{Y}$, the blowup of $Y$ along the base scheme $X$.

In \cite{Alu03}, Aluffi gave a procedure for computing the Segre class $i_* (X,\bbP^N)$ for a
closed subscheme $X \subset \bbP^N$.
The procedure is founded on the observation that the pushforward $\segre(X,\bbP^N)$
is determined by the class of the \emph{shadow} of the blowup of $\bbP^N$ along $X$.
This class is $G = a_0 [\bbP^0] + \dots + a_N [\bbP^N]$ where
\[a_i = \int \pi_* (c_1(\tilde{\mathrm{pr}}^* \calO(1))^{N-i} \capp [\tilde{Y}])\]
are the projective degrees of the map $pr$.
Once these coefficients have been found, the Segre class can be computed as (\cite[Theorem 2.1]{Alu03})
\[ i_* \segre(X,\bbP^N) = 1 - c(\calO(d))^{-1} \capp (G \otimes \calO(d)). \]
Helmer \cite{Hel14} has recently expanded on this approach by giving a new way to
compute the projective degrees $a_i$.

In \cite{EJP13}, Eklund, Jost, and Peterson gave a different procedure (call it EJP) for computing
$i_* \segre(X,\bbP^N)$ based on computing a residual intersection.
If $g_1, \dots, g_{N-r}$ are polynomials of degree $d$ which vanish on $X$,
then we can consider the linear system generated by the corresponding sections
in $H^0(\calI_X \otimes \calO(d), \bbP^N)$.
The base locus of this linear system contains $X$ and so we may write the intersection
product of the corresponding hypersurfaces as a part contributed by $X$
along with a ``residual'' part:
\[ D_1 \cdot \ldots \cdot D_{N-r} = \{c(\mathcal{N}) \capp \segre(X,\bbP^N)\}_{N-r} + \mathbf{R}_0 \]
The method of EJP is to sequentially add polynomials to get $g_1,\dots,g_m$
for $N-r \leq m \leq N$, getting residual contributions $\mathbf{R}_0, \dots, \mathbf{R}_r$.
They then use the formula \cite[Theorem 3.2]{EJP13}
\[ \sum_{i=0}^p \binom{m}{p-i} d^{p-i} s_{r-i} = d^m - \deg(\mathbf{R}_p),  \]
where $p = m-(N-r)$, to write the coefficients of $i_* \segre(X,\bbP^N)$ in terms of the
degrees of these residual intersection classes.

From the equations of EJP one can construct the following linear system for the $s_i$ in terms of $N,d,r$ and the residuals:
\begin{equation}
  \arraycolsep=3pt
  \medmuskip=1mu
  \label{eqn:EJPlinearsystem}
  \begin{pmatrix}
    1 & \binom{N}{1} d & \binom{N}{2} d^2 & \dots & \binom{N}{r} d^r \\[1ex]
    0 & 1 & \binom{N-1}{1} d & \dots & \binom{N-1}{r-1} d^{r-1} \\[1ex]
    \vdots & \vdots & \ddots & \ddots & \vdots \\[1ex]
    0 & 0 & \dots & 1 & \binom{N-r+1}{1} d \\[1ex]
    0 & 0 & \dots & 0 & 1
  \end{pmatrix}
  \begin{pmatrix}
    s_0 \\[1ex] s_1 \\[1ex] \vdots \\[1ex] s_{r-1} \\[1ex] s_r
  \end{pmatrix} = \begin{pmatrix}
    d^{N} - \deg(\mathbf{R}_r) \\[1ex] d^{N-1} - \deg(\mathbf{R}_{r-1}) \\[1ex] \vdots \\[1ex] d^{N-r+1} - \deg(\mathbf{R}_1) \\[1ex] d^{N-r} - \deg(\mathbf{R}_0)
  \end{pmatrix}
\end{equation}

\Cref{alg:main} can be seen as a generalization of EJP in the sense that the system of equations
produced using \cref{alg:main} exactly matches \cref{eqn:EJPlinearsystem} when $Y = \bbP^N$.
Indeed, Helmer has observed \cite[Section 3.2]{Hel14} that the projective degrees of $pr$ actually
correspond with the residual degrees, that is,
$a_{N-r+j} = \deg(\mathbf{R}_{j})$ for $0 \leq j \leq r$, and performing this substitution
produces \cref{eqn:almostMainlinearsystem} on page \pageref{eqn:almostMainlinearsystem}.
However, the geometric pictures in the two algorithms are different, even in the case $Y = \bbP^N$, since EJP works by constructing a system of equations which constrain $s(X,Y)$ and \cref{alg:main} works by producing several Segre classes and comparing them.

\section{Segre classes commute with intersecting by effective Cartier divisors}
\label{sec:segreClassesCommuteWithDivisors}
The main goal of this section is \cref{thm:segreClassesCommuteWithDivisors}.
This result can be found in many places, for instance \cite[Lemma 4.1]{Alu14},
but a detailed proof was harder to find so we give one here.

For this section, we will first need to review the concept of
distinguished varieties.
Given a closed embedding $X \subset Y$, let $C_X Y$ denote the normal cone to $X$ in $Y$.
Let $C_1, \dots, C_l$ be the irreducible components of $C_X Y$ and let $Z_i$
denote the support of $C_i$ in $X$.
These $Z_i$ are called the \emph{distinguished varieties} for $X$ in $Y$.
It should be noted that the $Z_i$ are not necessarily distinct; two irreducible components
$C_i, C_j$ may have the same support $Z_i = Z_j$.

Now let $K \subset Y$ be a Cartier divisor
and let $\tilde{K}$ be the blowup of $K$ along $X \cap K$, with exceptional divisor $F$.
Then $F = E \cap \tilde{K}$ and we have the following commutative diagram,
where the front and back faces are blowup squares.

\begin{minipage}{\linewidth}
  \centering
  \begin{tikzcd}[ampersand replacement=\&]
    \& E \arrow{rr} \arrow{dd}[pos=0.25, font=\small]{\pi^{\prime}} \& \& \tilde{Y} \ar{dd}[pos=0.25,font=\small]{\pi} \\
    F
    \ar{dd}[pos=0.25,font=\small]{\rho^{\prime}}
    \ar[crossing over]{rr} \ar[font=\small]{ur}
    \& \& \tilde{K} \ar{ur} \\
    \& X \ar[pos=0.25,font=\small]{rr} \& \& Y \\
    X \cap K \ar[pos=0.25]{rr} \ar[font=\small]{ur} \& \& K \ar{ur}
    \ar{uu}[crossing over, leftarrow,font=\small, swap, pos=0.25]{\rho}
  \end{tikzcd}
\end{minipage}

\begin{prop}
  \label{prop:distvars}
  Let $K \subset Y$ be an effective Cartier divisor.
  If $K$ does not contain any of the distinguished varieties for $X \subset Y$,
  then the proper transform $\tilde{K}$ is equal to the total transform $\pi^{-1}(K)$.
  \begin{proof}
    It suffices to show that each irreducible component of $\pi^{-1}(K)$ is contained in $\tilde{K}$.
    Aside from those in $\tilde{K}$, the irreducible components of $\pi^{-1}(K)$ each lie in $\pi^{-1}(K \cap Z)$ for some distinguished variety $Z$. But $K \cap Z$ is an effective Cartier divisor on $Z$, so its pullback is one on $\pi^{-1}(Z)$ and must have codimension at least 2 in $\tilde{Y}$. Thus since $\pi^{-1}(K)$ is pure of codimension 1 in $\tilde{Y}$, there is no component of $\pi^{-1}(K)$ contained in $\pi^{-1}(K \cap Z)$.
  \end{proof}
\end{prop}

We can now prove the formula below which shows how the Segre classes $\segre(X,Y)$
and $\segre(X\cap K, K)$ are related for divisors $K \subset Y$ which do not contain the distinguished
varieties of $X$.

\begin{thm}
  \label{thm:segreClassesCommuteWithDivisors}
  For an effective Cartier divisor $K \subset Y$ which does not contain any of the distinguished varieties for $X$ in $Y$,
  \[ K \cdot \segre(X,Y) = \segre(X \cap K, K) \in A_*(X \cap K). \]
  \begin{proof}
    By definition, we have 
    $s(X,Y) = \pi'_* \left( \frac{E}{1+E} \frown [\tilde{Y}] \right)$ and $s(X \cap K, K) = \rho'_* \left( \frac{E}{1 + E} \frown [\tilde{K}] \right)$.
    We arrive at the desired equation by combining \cref{prop:distvars}
    with the projection formula for divisors and commutativity of Chern classes:
    \begin{align*}
      K \cdot s(X,Y) &= K \cdot \pi'_* \left( \frac{E}{1+E} \frown [\tilde{Y}] \right) \\
      &= \rho'_* \left( \pi'^*(K) \left( \frac{E}{1+E} \frown [\tilde{Y}] \right) \right) \\
      &= \rho'_* \left( \left( \frac{E}{1+E} \frown [\tilde{K}] \right) \right) \\
      &= s(X \cap K, K). \qedhere
    \end{align*}
  \end{proof}
\end{thm}

We now give the precise statement that will be needed for \cref{alg:main}.

\begin{corollary}
  \label{cor:segreClassesCommuteWithDivisors}
  For a general hypersurface $H \subset \bbP^N$ of any fixed degree,
  \[ H \cdot \segre(X,Y) = \segre(X \cap H, Y \cap H) \in A_*(X \cap H). \]
  \begin{proof}
    Bertini's Theorem says that for general $H \subset \bbP^N$, the divisor $H \cap Y$ is effective and does not contain any distinguished varieties for $X \subset Y$, so \cref{thm:segreClassesCommuteWithDivisors} applies.
  \end{proof}
\end{corollary}

\section{Degrees of Segre classes via linear projection}
\label{sec:segreDegrees}

Now assume $X$ is the base scheme of a ($k+1$ dimensional) linear system
$L \subset H^0(Y, \calO(d))$ and let $\mathrm{pr}: Y \dashrightarrow \bbP^k$ be the associated rational map.
As usual, denote by $\pi: \tilde{Y} \rightarrow Y$ the blowup of $Y$ along $X$ and let
$\tilde{\mathrm{pr}}: \tilde{Y} \rightarrow \bbP^k$ be the lift of the projection map to $\tilde{Y}$.
In this setup, \cite[Proposition 4.4]{Ful98} gives the formula
\[ \int\limits_X c(\calO(d))^n \capp \segre(X,Y) = \int\limits_Y c_1(\calO(d))^n - \deg(\tilde{Y} / \tilde{\mathrm{pr}}(\tilde{Y})) \int\limits_{\bbP^k} c_1(\calO(1))^n \capp [\tilde{\mathrm{pr}}(\tilde{Y})]. \]
If the pushforward of $s(X,Y)$ to $\bbP^N$ is $i_* \segre(X,Y)=s_0 [\bbP^0] + s_1 [\bbP^1] + \dots + s_r [\bbP^r]$
then the left-most degree is
\begin{equation}
  \label{eqn:segreintegralexpansion}
  \left\{ (1 + d[\bbP^{N-1}])^n (s_0 [\bbP^0] + s_1 [\bbP^1] + \dots + s_r [\bbP^r]) \right\}_0
  = \sum_{i=0}^r \binom{n}{i} d^i s_i
\end{equation}
while for the middle term we have $\int_Y c_1(\calO(d))^n = d^n \deg(Y)$.

As for the last term, we can interpret the factors as follows.
The degree $\deg(\tilde{Y} / \tilde{\mathrm{pr}}(\tilde{Y}))$ is the degree of the map
$\tilde{\mathrm{pr}}$ when this is generically finite, and is 0 otherwise,
while $\int_{\bbP^k} c_1(\calO(1))^n \capp [\tilde{\mathrm{pr}}(\tilde{Y})]$
is the degree of $\tilde{\mathrm{pr}}(\tilde{Y}) \subset \bbP^k$ if this subvariety
has dimension $n = \dim Y$ and is 0 otherwise.
Since the second term will yield 0 exactly when the first does, we can consolidate these by writing
\[ \deg(\tilde{Y} / \tilde{\mathrm{pr}}(\tilde{Y})) \cdot \int_{\bbP^k} c_1(\calO(1))^n \capp [\tilde{\mathrm{pr}}(\tilde{Y})]
  = \deg(\tilde{\mathrm{pr}}) \deg(\tilde{\mathrm{pr}}(\tilde{Y})) \]
with the understanding that the $\deg(\tilde{\mathrm{pr}})$ should be 0 if the map is not generically finite.

Altogether, this lets us write
\[ \sum_{i=0}^{\dim{X}} \binom{n}{i} d^i s_i
  = d^n \deg(Y) - \deg(\tilde{\mathrm{pr}}) \deg(\tilde{\mathrm{pr}}(\tilde{Y})) \]
Since the map $\tilde{\mathrm{pr}}$ factors as $\mathrm{pr} \circ \pi$ and $\pi$ is an
isomorphism away from $X$, we see that $\deg(\tilde{\mathrm{pr}}) = \deg(\mathrm{pr})$
and $\tilde{\mathrm{pr}}(\tilde{Y}) = \overline{\tilde{\mathrm{pr}}(\tilde{Y} - E)}
= \overline{\mathrm{pr}(Y - X)}$.
This allows us to completely remove consideration of $\tilde{Y}$ in the computation of
$\segre(X,Y)$.

We collect these observations in the following proposition.

\begin{prop}
  \label{prop:segredegreeformula}
  In the setup of this section,
  \[ \sum_{i=0}^{\dim{X}} \binom{n}{i} d^i s_i
    = d^n \deg(Y) - \deg(\mathrm{pr}) \deg(\overline{\mathrm{pr}(Y-X)}).
    \xqedhere{1.8cm}{\qed} \]
\end{prop}

This should be seen as a linear equation in the unknowns $s_0,\dots,s_{\dim{X}}$.

\section{An algorithm for computing Segre classes}
\label{sec:algorithm}
In this section, we finally describe the algorithm for computing the Segre class
${i_* s(X,Y) \in A_*(\bbP^N)}$.

\subsection{Setup}
Start by setting $X_0 = X$ and $Y_0 = Y$
and assume that $X_0$ is cut out in $Y_0$ by homogeneous equations $f_1,\dots,f_m$.
If $d = max_i (\deg(f_i))$, we can consider $X_0$ as the base locus of a linear system in
$H^0(Y_0,\calO(d))$ by ``padding'' the $f_i$ appropriately.
That is, if $\deg(f_i) < d$ we replace $f_i$ by the polynomials $g_{i,1}, \dots, g_{i,q}$
where $g_{i,j} = f_i \sigma_j$ and $\sigma_1, \dots, \sigma_q$ are
the standard monomials of degree $d - \deg(f_i)$ in the coordinates of $\bbP^N$.
Relabeling, write the set of $g_{i,j}$ for all $i,j$ as $g_0, \dots, g_k$.
Then $X_0$ is the base locus of the linear system generated by the sections
corresponding to these $g_i$.

\subsection{Linear system of equations}
To shorten the statements that follow, define
$\delta_j = d^{\dim{Y_j}} \deg(Y_j) - \deg(\mathrm{pr}_j) \deg(\overline{\mathrm{pr_j}(Y_j-X_j)})$,
where $\mathrm{pr}_0: Y_0 \dashrightarrow \bbP^k$ is the projection defined previously and $\mathrm{pr}_j$ will be defined analogously for objects $X_j \subseteq Y_j$.
Writing $i_* \segre(X_0,Y_0) = s_0 [\bbP^0] + s_1 [\bbP^1] + \dots + s_r [\bbP^r]$ as before,
we get \[ \sum_{i=0}^r \binom{n}{i} d^i s_i  = \delta_0 \] by \cref{prop:segredegreeformula}.

Now choose a general hyperplane section $H \subset Y_0$ \`{a} la \cref{cor:segreClassesCommuteWithDivisors}.
If we set $Y_1 = Y_0 \cap H$ and $X_1 = X_0 \cap H$ we can repeat the procedure above,
with the sections of the new linear system given by restricting $g_0, \dots, g_k$
to $Y_1$.
This gives a new projection $pr_1$.
By \cref{cor:segreClassesCommuteWithDivisors}, we have $H \cdot s(X_0,Y_0) = s(X_1,Y_1)$
so the resulting linear equation will be
\[ \sum_{i=1}^{r} \binom{n-1}{i-1} d^{i-1} s_i = \delta_1. \]

Repeating this process of taking hyperplane sections will yield a system of equations
$A \vec{s} = \vec{\delta}$.  Explicitly,
\begin{equation}
  \label{eqn:almostMainlinearsystem}
  \begin{pmatrix}
    1 & \binom{n}{1}d & \binom{n}{2} d^2 & \dots & \binom{n}{r} d^r \\[1ex]
    0 & 1 & \binom{n-1}{1}d & \dots & \binom{n-1}{r-1} d^{r-1} \\[1ex]
    \vdots & \vdots & \ddots & \ddots & \vdots \\[1ex]
    0 & 0 & \dots & 1 & \binom{n-r+1}{1} d \\[1ex]
    0 & 0 & \dots & 0 & 1
  \end{pmatrix}
  \begin{pmatrix}
    s_0 \\[1ex] s_1 \\[1ex] \vdots \\[1ex] s_{r-1} \\[1ex] s_r
  \end{pmatrix} = \begin{pmatrix}
    \delta_0 \\[1ex] \delta_1 \\[1ex] \vdots \\[1ex] \delta_{r-1} \\[1ex] \delta_r
  \end{pmatrix}
\end{equation}
and solving for $A^{-1}$ gives
\begin{equation}
  \arraycolsep=3pt
  \medmuskip=1mu
  \label{eqn:mainlinearsystem}
  \begin{pmatrix}
    1 & \binom{N}{1} (-d) & \binom{N}{2} (-d)^2 & \dots & \binom{N}{r} (-d)^r \\[1ex]
    0 & 1 & \binom{N-1}{1} (-d) & \dots & \binom{N-1}{r-1} (-d)^{r-1} \\[1ex]
    \vdots & \vdots & \ddots & \ddots & \vdots \\[1ex]
    0 & 0 & \dots & 1 & \binom{N-r+1}{1} (-d) \\[1ex]
    0 & 0 & \dots & 0 & 1
  \end{pmatrix}
  \begin{pmatrix}
    \delta_0 \\[1ex] \delta_1 \\[1ex] \vdots \\[1ex] \delta_{r-1} \\[1ex] \delta_r
  \end{pmatrix} = \begin{pmatrix}
    s_0 \\[1ex] s_1 \\[1ex] \vdots \\[1ex] s_{r-1} \\[1ex] s_r
  \end{pmatrix}.
\end{equation}

From this we get the following closed form expression for the $s_i$.

\begin{prop}
  \label{prop:segreformula}
  With notation as in this section,
  \begin{equation*}
    s_i = \sum_{j=0}^r \binom{N-i}{j-i} (-d)^{j-i} \delta_j.
    \xqedhere{4.11cm}{\qed}
  \end{equation*}
\end{prop}

\subsection{Algorithm}
The algorithm can now be summarized quite concisely:
\begin{enumerate}[label={Step \arabic*:}]
  \item Calculate $\deg(Y), \deg(pr)$ and $\deg(\overline{pr(Y-X)})$.
  \item Replace $X$ and $Y$ with $X \cap H$ and $Y \cap H$ for a suitable hyperplane $H \subset \bbP^N$ and repeat Step 1.  Continue until $\dim X = 0$.
  \item For each $i$, calculate $s_i$ using \cref{prop:segreformula}.
\end{enumerate}

As a more verbose presentation we offer the pseudocode version of the algorithm for computing $\segre(X,Y)$.
In practice, the inputs are two ideals, $I_X, I_Y \subset k[x_0,\dots,x_N]$
defining the schemes $X,Y \subset \bbP^N_k$ respectively.

\begin{algorithm}[h]
  \caption{Computation of s(X,Y)}
  \label{alg:main}
  \DontPrintSemicolon
  \KwIn{a variety $Y^n \subset \mathbb{P}^N$ and closed embedding $X^r \rightarrow Y$ given by a linear system in $\mathcal{O}(d)$}
  \KwOut{the Segre class $\segre(X,Y)$ pushed forward to $\mathbb{P}^N$}
  \BlankLine
  $X_0 \gets X$, \: $Y_0 \gets Y$\;
  \For{$i=0$ \KwTo $r$}{
    $\delta_i \gets d^n \operatorname{deg}(Y_i) -
    \operatorname{deg}(\mathrm{pr_i})
    \operatorname{deg}(\overline{\mathrm{pr_i}(Y_i - X_i)}) $\;
    Choose generic hyperplane $H \subset \mathbb{P}^N$ as in \cref{cor:segreClassesCommuteWithDivisors}\;
    $Y_{i+1} \gets Y_i \cap H$, \: $X_{i+1} \gets X_i \cap H$\;
  }
  \For{$i=0$ \KwTo $r$}{
    $s_i \gets \sum_{j=0}^r \binom{N-i}{j-i} (-d)^{j-i} \delta_j$\;
  }
  \KwRet{$s_0 \mathbf{h}^N + \dots + s_r \mathbf{h}^{N-r}$}\;
\end{algorithm}
The algorithm has been implemented in Macaulay2 and an implementation should soon be available for SageMath as well.\footnote{\url{http://github.com/coreysharris/FMPIntersectionTheory-M2}}
The Macaulay2 package interfaces with Schubert2, as can be seen in the examples to follow.

\section{Applications}
\label{sec:Applications}

In this section we give examples of the general types of computations that can be made
if Segre classes are accessible.
Along with Chern classes (like that of Mather, Schwartz-MacPherson) one can also compute
polar classes, which can be used to compute Euclidean distance degrees.
We also discuss how Segre classes can be applied in the computation of degrees of projections of degeneracy loci.
Finally, we give several examples of intersection product computations, the flagship application of the Segre class.

\subsection{Chern classes}
Let $Z$ be a hypersurface of a nonsingular $m$-dimensional variety $M$ and let $J$ be the singularity subscheme of $Z$, i.e., the subscheme defined by the Jacobian ideal of $Z$.
Several important objects instrinsic to $Z$ can be computed if the Segre class of $J$ in $Z$ or $M$ is known.

In \cite[Lemma I.4]{Alu99}, Aluffi gave a formula for the Chern-Schwartz-MacPherson class of $Z$ in terms of $\segre(J,M)$:
\[ c_{SM}(Z) = c(TM) \capp \left( \frac{[Z]}{1+Z} + (c(\calO(Z)|_J) \capp \segre(J,M))^\vee \otimes \calO(Z)|_J \right). \]
Similarly, the Chern-Mather class of $Z$ can be written (see \cite[Prop. 2.2]{AB03}	)
\[ c_{Ma}(Z) = c(TM) \capp \left( \frac{[Z]}{1+Z} + \segre(J,Z)^\vee \otimes \calO(Z)|_J \right). \]
where now the key input is $\segre(J,Z)$.

Both of these formulas are quite explicit, and can be implemented in Macaulay2 \cite{M2}
by a literal translation of the formulas above.
For instance, Schubert2 provides the structure sheaf {\ttfamily O} and also provides {\ttfamily adams} which performs the $(\cdot)^\vee$ operation,
and our code implements $\otimes$ as {\ttfamily **}.
Then if {\ttfamily cZ} is the cycle class of $Z$ in $A_*(\bbP^N)$, the Chern-Mather class can be computed with

\begin{minipage}{0.95\linewidth}
\begin{lstlisting}
s = segreClass(iJ,iZ)
chern(TM)*(cZ*(1+cZ)^(-1)+(adams(-1,s)**O(cZ)))
\end{lstlisting}
\end{minipage}

For example, we can compute $c_{SM}(X)$ and $c_{Ma}(X)$ for $X \subset \bbP^3$ the singular cone:

\begin{minipage}{0.95\linewidth}
\begin{lstlisting}
i1 : PP3 = QQ[x,y,z,t]
i2 : X = ideal x*y-z^2
i3 : chernMather(X)
o3 = 2H    + 4H    + 2H
       2,1     2,2     2,3
i4 : chernSchwartzMacPherson(X)
o4 = 2H    + 4H    + 3H
       2,1     2,2     2,3
\end{lstlisting}
\end{minipage}
\\The output of these two commands are elements of the {\ttfamily IntersectionRing} of $\bbP^3$.
In this case they say $c_{SM}(X) = 2[\bbP^2] + 4[\bbP^1] + 3[\bbP^0]$ and
$c_{Ma}(X) = 2[\bbP^2] + 4[\bbP^1] + 2[\bbP^0]$.

\subsection{Polar classes}

Now assume $Z$ is a hypersurface in $\bbP^N$ and
let $\Lambda_{k-1} \subset \bbP^N$ be a general $(k-1)$-dimensional linear subspace.
We define the $k$th \emph{polar locus} of $Z$ (relative to $\Lambda_{k-1}$) by
\[P_k = \overline{\{ x \in Z_{sm} \;|\; \Lambda_{k-1} \subset T_x X\}}.\]
The $k$th \emph{polar class} $\varrho_k$ is the degree $\varrho_k = \int [P_k]$.
Piene \cite[Theorem 2.3]{Pie78} showed that when $Z$ is a hypersurface of degree $d$, the class $[P_k] \in A_{N-1-k}(Z)$ can be expressed
in terms of the Segre class of the singularity subscheme ${s(J,Z) = s_0[\bbP^0] + \dots + s_{N-1}[\bbP^{N-1}]}$:
\[
  [P_k] = (d-1)^k h^k \capp [Z] - \sum_{i=0}^{k-1} \binom{k}{i} (d-1)^i h^i \capp s_{i}[\bbP^{i}] \in A_{N-1-k}(\bbP^N)
\]
The polar classes $\varrho_0,\dots,\varrho_{N-1}$ are the degrees of the classes $[P_k]$, so when $Z$ is a hypersurface the above formula gives

\begin{equation}
  \label{eqn:polarclasses}
  \varrho_k = \int [P_k] = d(d-1)^k - \sum_{i=0}^{k-1} \binom{k}{i} (d-1)^i s_{i}.
\end{equation}

Polar classes are projective invariants \cite{Pie78,Hol88}.
If $Z$ has codimension $c > 1$ in $\bbP^N$ and $pr_\Lambda: \bbP^N \dashrightarrow \bbP^{N-c+1}$
is the projection from a general $(c-1)$-dimensional linear subspace $\Lambda \subset \bbP^N$
then $\varrho_k(Z) = \varrho_k(pr_\Lambda(Z))$ for all $k$.
Thus \cref{alg:main} allows the computation of the polar classes of any projective variety.

Moreover, Piene \cite[p.~19]{Pie88} showed that the polar classes carry the same information
as the Chern-Mather class.
In particular, if $\varsigma_k$ is the coefficient of $[\bbP^k]$ in $c_{Ma}(Z)$,
\[ \varsigma_k = \sum_{i=0}^k \binom{n+1-i}{k-i} \varrho_k \quad \text{ and } \quad \varrho_k = \sum_{i=0}^k \binom{n+1-i}{k-i} \varsigma_k. \]

An implication of this is that Chern-Mather classes are also projective invariants, and indeed
Piene observed that if $pr$ is a general linear projection as above,
\[c_{Ma}(pr(Z)) = {pr}_* c_{Ma}(Z).\]
in $A_*(pr(Z)).$

\subsection{Euclidean distance degrees}
For any subvariety $Z \subset M$,
an alternate but equivalent characterization of its polar classes is
the multidegree of the conormal cycle of $Z$ in $M$.
That is, let
\[ \mathcal{N}_Z^\circ = \left\{ (p,h) \in \bbP^N \times (\bbP^N)^\vee \;|\; p \in Z_{sm} \text{ and $h$ is tangent to $X$ at $p$}  \right\}. \]
The closure $\mathcal{N}_Z = \overline{\mathcal{N}_Z^\circ} \subset \bbP^N \times (\bbP^N)^\vee$ is the \emph{conormal variety} of $Z$.
It is irreducible of pure dimension $N-1$ and its class in the Chow group $A_*(\bbP^N \times (\bbP^N)^\vee)$ can be written
\[ [\mathcal{N}_Z] = \varrho_0 a^N b + \varrho_1 a^{N-1} b^2 + \dots + \varrho_{N-1} a b^{N}  \]
where $a,b$ are the pullbacks of the hyperplane classes on $\bbP^N$ and $(\bbP^N)^\vee$ respectively.

For a general point $y \in \bbP^N$ define $d_y(x) = \sum_{i=0}^N (y_i-x_i)^2$ to be the
squared Euclidean distance from $y$ to $x$.
The \emph{Euclidean distance degree} $\mathrm{EDdeg}(Z)$ of an algebraic variety $Z \subset \bbP^N$
is the number of critical points for $d_y$ restricted to $Z_{sm}$.

In \cite[Theorem 5.4]{DHO+15}, the authors show that if $\mathcal{N}_Z$ does not meet the diagonal in $\bbP^N \times (\bbP^N)^\vee$, then $\mathrm{EDdeg}(Z)$ is the sum of the polar classes of $Z$.

For example (\hspace{1sp}\cite[Example 5.7]{DHO+15}), consider Cayley's quartic surface in $C \subset \bbP^3$ defined by $x^3-xy^2-xz^2+2yzt-xt^2 =0$.
The conormal variety is defined by 18 equations in $\bbP^3 \times \bbP^3$, whereas the Segre class of the singularity subscheme is quite simple:

\begin{minipage}{0.95\linewidth}
\begin{lstlisting}
i1 : C = ideal "x3-xy2-xz2+2yzt-xt2"
i2 : J = ideal singularLocus C
i3 : segreClass(J,C)
o3 = 8H
       2,3
\end{lstlisting}
\end{minipage}
\\ so $\segre(J,C) = 8[\bbP^0]$ and thus by \cref{eqn:polarclasses}, the Euclidean distance degree is $3+6+4 = 13$ (cf.~\cite{DHO+15}).

\subsection{Degrees of projections of degeneracy loci}

Let $\bbP^{n^2-1}$ be the space of $n \times n$ complex matrices.
Let $\tau_k$ denote the locus of matrices of corank at least $k$.
In \cite{Alu14},
Aluffi considers the following problem: if $S$ is a set of $s$ coordinates on $\bbP^{n^2-1}$,
and $\pi_S$ is the projection from center $V(S)$,
what is the degree of the closure of $\pi_S(\tau_k)$ in $\bbP^{n^2-s-1}$?
An excess intersection computation yields the general solution (cf. \cite[Proposition 4.4]{Ful98}).
\begin{lemma}[\hspace{1sp}{\cite[Lemma 2.7]{Alu14}}]
  Let $H$ be the class of a hyperplane in $\bbP^{n^2-1}$.
  Then the degree $d$ of the closure of $\pi_S(\tau_k)$ is given by
  \[ d = \prod_{i=0}^{k-1} \frac{\binom{n+i}{k}}{\binom{k+i}{k}} - \int (1+H)^{n^2-k^2-1} \capp \segre(L_S \cap \tau_k, \tau_k). \xqedhere{2.35cm}{\qed} \]
\end{lemma}
These Segre classes are great examples of problems for which a general computational method
was not previously available.
As an example of the efficacy of \cref{alg:main}, the table below shows some computations
for $\tau_2 \subset \bbP^8$, the locus of $3 \times 3$ matrices of rank $\leq 2$.
The tableaux denote the sets $S$, where the top left square is the
coordinate of the first row and first column of the general matrix.

\newcommand\Tstrut{\rule{0pt}{3.4ex}}         
\newcommand\Bstrut{\rule[-2.8ex]{0pt}{0pt}}   
\newlength{\plen}
\newlength{\mlen}
\settowidth{\plen}{${}+{}$}
\settowidth{\mlen}{${}-{}$}
\newcommand{\mc}[2]{\multicolumn{#1}{r@{\hspace{\plen}}}{#2}}
\newcommand{\mcminus}[1]{\multicolumn{1}{r@{{}-{}}}{#1}}

\[
  \begin{array}{|c|r@{{}+{}}r@{{}+{}}r@{{}+{}}r@{{}+{}}r|}
    \hline
    S & \multicolumn{5}{|c|}{i_*\segre(L_S \cap \tau_2, \tau_2)}\\
    \hline
    {\tiny\Yvcentermath1 \yng(3,1)} & -25[\bbP^0] & \mcminus{10[\bbP^1]} & \mcminus{2[\bbP^2]} & [\bbP^3] & [\bbP^4] \Tstrut\\[0.5em]
    {\tiny\Yvcentermath1 \yng(3,2)} & \mc{1}{10[\bbP^0]} &  \mcminus{} & 2[\bbP^2] & \mc{1}{[\bbP^3]} & \\[0.5em]
    {\tiny\Yvcentermath1 \yng(3,2,1)} & \mcminus{3[\bbP^0]} & 3[\bbP^1] & \mc{1}{[\bbP^2]} & \mc{1}{} & \Bstrut \\[0.5em]
    {\tiny\Yvcentermath1 \young(\hfil\hfil\hfil,\hfil\hfil,::\hfil) } & \mcminus{3[\bbP^0]} & 3[\bbP^1] & \mc{1}{[\bbP^2]} & \mc{1}{} & \Bstrut \\[0.5em]
    {\tiny\Yvcentermath1 \yng(3,3,1)} & -11[\bbP^0] & \mc{1}{2[\bbP^1]} & \mc{2}{} & \Bstrut \\[0.5em]
    \hline
  \end{array}
\]

As an example we compute the first row in the table:

\begin{minipage}{0.95\linewidth}
\begin{lstlisting}
i1 : PP8 = QQ[v_0..v_8]
i2 : I = minors(3, matrix{{v_0..v_2},{v_3..v_5},{v_6..v_8}})
i3 : J = ideal (v_0..v_3)
i4 : segreClass(J,I)  -- J is replaced by J+I automatically
o4 = H    - H    - 2H    + 10H    - 25H
      2,4    2,5     2,6      2,7      2,8
\end{lstlisting}
\end{minipage}

\subsection{Intersection Products}
\label{subsec:intersectionProducts}

A basic intersection product setup in \cite[Ch. 6]{Ful98} is a fiber square \\
\begin{minipage}{\linewidth}
  \centering
  \begin{tikzcd}[ampersand replacement=\&]
    \& W \ar{r}[font=\small]{j} \ar{d}[pos=0.5, font=\small]{g} \& V \ar{d}[pos=0.5,font=\small]{f} \\
    \& X \ar{r}[font=\small]{i} \& Y
  \end{tikzcd}
\end{minipage}
in which $i$ is a regular embedding and $V$ is a variety of dimension $k$.
In this situation, we let $N$ denote the pullback $g^* N_X Y$ of the normal bundle to $X$ in $Y$
and define the intersection product of $X$ and $V$ in $Y$ to be
\[ X \cdot_Y V = \{ c(N) \capp \segre(W,V) \}_{k+r-n}. \]
Assuming that $f,i$ are really inclusions with $X,V \subset Y \subset \bbP^N$,
we potentially have everything needed to compute pushforward of $X \cdot Y$ to $\bbP^N$.
That is, if $W \overset{\eta}\hookrightarrow \bbP^N$ is the inclusion and we can find a bundle $\hat{N}$ on $\bbP^N$ such that $N = \eta^* \hat{N}$, then we can use the projection formula
\begin{align*}
  \eta_*(X \cdot V) &= \eta_* \left(\{ c(N) \capp \segre(W,V) \}_{k+r-n} \right) \\
  &= \{ c(\hat{N}) \capp \eta_* \segre(W,V) \}_{k+r-n}.
\end{align*}

Such is the situation when $X \subset \bbP^N$ and $Y \subset \bbP^N$ are both complete intersections.  Then the normal bundle $N_X Y$ is isomorphic to the quotient $N_X \bbP^N / i^* N_Y \bbP^N$ and so its total chern class in $\bbP^N$ can be written
\[\frac{\prod_i (1+d_i [\bbP^{4}])}{\prod_j (1+e_j [\bbP^{4}])}.\]

\subsubsection{Grassmannian of lines in $\bbP^3$}
As an example, consider $G(2,4) \subset \bbP^5$, the Grassmannian of lines in $\bbP^3$,
defined by the equation $ab-cd+ef=0$.
Let $\Sigma_{1}$ be the Schubert cycle of lines meeting a fixed line, 
with representative defined by the equation $b=0$ on $G(2,4)$, 
and let $\Sigma_{2,1}$ be the Schubert cycle of lines containing a fixed point and lying in a fixed plane, 
with representative defined by the equations $b=d=e=f=0$.
Say we wish to compute the pushforward to $\bbP^N$ of the intersection product $\Sigma_{2,1} \cdot \Sigma_1$ in $G(2,4)$.
By the remarks above, we can write the chern class of the normal bundle as
\[ N_{\Sigma_{2,1}} G(2,4) = \frac{(1+[\bbP^{4}])^4}{(1+2[\bbP^{4}])} \]
and the other needed ingredient is $\segre(\Sigma_1 \cap \Sigma_{2,1}, \Sigma_1) = \segre(\Sigma_{2,1}, \Sigma_1)$.  In Macaulay2, we compute this by

\begin{minipage}{0.95\linewidth}
\begin{lstlisting}
i1 : PP5 = QQ[a,b,c,d,e,f]
i2 : G = ideal "ab-cd+ef"
i3 : S1 = G + ideal b
i4 : S21 = G + ideal (b,c,d,e)
i5 : segreClass(S21,S1)
o5 =  H    - H
       2,4    2,5
\end{lstlisting}
\end{minipage}
which returns $\segre(\Sigma_{2,1},\Sigma_1) = -[\bbP^0] + [\bbP^1]$.
Thus, the (pushed-forward) intersection product is
\begin{align*}
  \Sigma_{2,1} \cdot \Sigma_1 &= \left\{ \frac{(1+[\bbP^{4}])^4}{(1+2[\bbP^{4}])} (-[\bbP^0] + [\bbP^1]) \right\}_0 \\
  &= \left\{ \left( -2[\bbP^0] + [\bbP^1] +  2[\bbP^3] + 2[\bbP^4] + [\bbP^5]  \right)
    \left( -[\bbP^0] + [\bbP^1] \right) \right\}_0 \\
  &= \left\{ [\bbP^0] + [\bbP^1] \right\}_0 \\
  &= [\bbP^0]
\end{align*}
as expected since $\sigma_{2,1} \sigma_{1} = \sigma_{2,2}$ in $A_*(G(2,4))$.

Now let $V$ be the subvariety of $G(2,4)$ defined by $b^2-cf=0$.
Then $\dim V = 3$ and $\dim V_{sing} = 2$.
If we wish to compute the pushforward of $\Sigma_1 \cdot V$ we can use
\begin{align*}
  \left\{ c(N_{\Sigma_1} G(2,4)) \capp \segre(\Sigma_1 \cap V, V) \right\}_2 &= \left\{ \frac{(1+[\bbP^{4}])(1+2[\bbP^{4}])}{1+2[\bbP^{4}]} \capp \segre(\Sigma_1 \cap V, V) \right\}_2 \\
  &= \left\{ (1+[\bbP^4])(4[\bbP^0] - 4[\bbP^1] + 4[\bbP^2]) \right\}_2 \\
  &= 4[\bbP^2].
\end{align*}
But this example is too easy, in the sense that it could have been computed as a product in $\bbP^5$.
That is, since $\Sigma_1$ is a complete intersection in $G(2,4)$, we could write $\Sigma_1 \cdot V$ as $G(2,4) \cdot H \cdot V$ where $H$ is the hypersurface in $\bbP^5$ defined by $b=0$.
Then Bezout's Theorem says the intersection product has degree 4, so counting dimensions yields the result.

In light of this we consider $X \subset G(2,4)$ defined by $b=d=f=0$.  This is still a complete intersection in $\bbP^5$, but is not a complete intersection in $G(2,4)$, so the intersection product $X \cdot V$ cannot just be directly computed in $\bbP^5$.
Notice also that $X \subset V$, so the intersection is not proper.
The computation this time is
\begin{align*}
  X \cdot V &= \left\{ c(N_X G(2,4)) \capp \segre(X \cap V, V) \right\}_1 \\
  &= \left\{ \frac{(1+[\bbP^{4}])^3}{1+2[\bbP^{4}]} \capp (-2[\bbP^0] + [\bbP^1] + [\bbP^2]) \right\}_1 \\
  &= 2[\bbP^1]
\end{align*}
We can still verify this result by writing the product in Schubert cycles.
Since $V$ is a hypersurface in $G(2,4)$ of degree 4, it must be $[V] = 2 \sigma_{1} \in A_*(G(2,4))$.
With a little more work, one can check that $[X] = \sigma_{1,1}$, and so we find
\[ X \cdot V = 2 \sigma_1 \sigma_{1,1} = 2 \sigma_{2,1} \]
and $\Sigma_{2,1}$ is a line in $\bbP^5$.

\bibliography{segre}{}
\bibliographystyle{amsalphaeprint}

\end{document}